\documentclass[12pt]{amsart}
\usepackage{amsmath}
\usepackage{relsize}
\usepackage{amssymb}
\usepackage{amsthm}
\usepackage{mleftright}
\usepackage{amscd}
\usepackage{textcomp}
\usepackage{color}
\usepackage[latin1]{inputenc}
\usepackage{mathpazo}
\usepackage[scaled=.95]{helvet}
\usepackage{courier}
\usepackage[all]{xy}
\usepackage{tikz} 
\usetikzlibrary{babel} 
\usepackage[hypcap]{caption}

\usetikzlibrary{calc} 

\usetikzlibrary{through} 
\usetikzlibrary{patterns} 
\usetikzlibrary{intersections} 
\usetikzlibrary{matrix} 
\usetikzlibrary{cd} 

\swapnumbers
\headsep=1cm \numberwithin{equation}{section}
\usepackage[colorinlistoftodos]{todonotes}

\newtheorem{theorem}{Theorem}[section]
\newtheorem{lemma}[theorem]{Lemma}
\newtheorem{proposition}[theorem]{Proposition}
\newtheorem{corollary}[theorem]{Corollary}

\theoremstyle{definition}
\newtheorem{definition}[theorem]{Definition}
\newtheorem{remark}[theorem]{Remark}
\newtheorem{remark and definition}{Remark and Definition}[section]
\newtheorem{remark and notation}{Remark and Notation}[section]

\newtheorem{example}[theorem]{Example}

\newcommand\Tor{\operatorname{Tor}}
\newcommand\height{\operatorname{height}}

\newcommand\grade{\operatorname{grade}}

\newcommand\Ker{\operatorname{\Ker}}
\newcommand\coker{\operatorname{coker}}

\newcommand\pd{\operatorname{pd}}

\newcommand\rank{\operatorname{rank}}

\newcommand{\test}[2]{\mleft(\knds\genfrac..{0pt}{}{#1}{#2}\knds\mright)}
\newcommand{\knds}{\kern-\nulldelimiterspace}

\begin{document}

\title{On Betti numbers for symmetric powers of modules}

\date{}


\author{V. H. Jorge-P\'erez }
\author{J. A. Santos-Lima}

\address{Universidade de S{\~a}o Paulo -
ICMC
Caixa Postal 668, 13560-970, S{\~a}o Carlos-SP, Brazil}

\email{vhjperez@icmc.usp.br}
\email{seyalbert@gmail.com}


\date{}
\thanks{The first author was supported by grant 2019/21181-0, S\~ao Paulo Research Foundation (FAPESP)  and the second author was supported by CAPES-Brasil Code 001.}


\begin{abstract}  
In this paper, we investigate the homological properties of the symmetric powers of finitely generated modules on a Noetherian local ring \( R \), focusing on the computation of minimal free resolutions and Betti numbers associated with \( \mathcal{S}_j(M) \), the \( j \)-th symmetric power of \( M \). Motivated by results of Weyman, Tchernev, Avramov, and Molica-Restuccia on the acyclicity of the complexes \( \mathcal{S}_j\textbf{F}_{\bullet} \), we study under what conditions the Betti numbers of \( \mathcal{S}_j(M) \) can be determined from those of \( M \). We prove that for modules \( M \) with a minimal finite free resolution satisfying the \( (SW_j) \) condition, the complex \( \mathcal{S}_j \textbf{F}_{\bullet} \) yields a minimal free resolution of \( \mathcal{S}_j(M) \) with finite projective dimension. We provide explicit formulas for the Betti numbers of \( \mathcal{S}_j(M) \) in terms of those of \( M \), including applications to linear-type modules and powers of ideals, as well as establishing upper and lower bounds for these invariants. As a consequence, we obtain special cases of the Buchsbaum-Eisenbud-Horrocks conjecture.  
\end{abstract}
\maketitle
\section{Introduction}

Throughout this paper, $(R, \mathfrak{m}, k)$ denotes a  Noetherian local ring where $\mathfrak{m}$ is the maximal ideal and $k:=R/\mathfrak{m}$ is the residue field, and every \( R \)-module \( M \) is finitely generated over \( R \). For an \( R \)-module \( M \), we denote by \( \mathcal{S}_j(M) \) the \( j \)-th symmetric power of \( M \), or equivalently, the \( j \)-th graded component of the symmetric algebra \( \mathcal{S}_R(M) \).  

Our main goal in this paper is to examine the homological properties related to the symmetric powers of finitely generated modules. Specifically, we are interested in computing the minimal free resolutions and Betti numbers of these symmetric powers. This investigation is motivated by the works of Weyman \cite{Weyman}, Tchernev \cite{Tchernev}, Avramov \cite{Avramov}, and Molica and Restuccia \cite{Molica}, where the acyclicity of the complexes \( \mathcal{S}_j \textbf{F}_{\bullet} \) associated with \( \mathcal{S}_j(M) \) is studied. Due to the acyclicity criteria established in \cite{Weyman, Tchernev, Avramov, Molica}, we observe that not all modules of finite projective dimension admit symmetric powers whose minimal free resolutions arise directly from the minimal resolution of \( M \). This leads us to the following natural question. 

{\bf Question:} Is it possible to determine the Betti numbers of \( \mathcal{S}_j(M) \) from the Betti numbers of \( M \)?  

We now outline the structure of the paper. Sections \ref{section2} and \ref{section3} present foundational results on symmetric algebras, symmetric powers, and divided powers.  

In Section \ref{section4}, we study the minimality of the complex \( \mathcal{S}_j \textbf{F}_{\bullet} \). For instance, we show that if \( M \) has a minimal finite free resolution and satisfies the \( (SW_j) \) condition, then \( \mathcal{S}_j \textbf{F}_{\bullet} \) is a minimal free resolution of \( \mathcal{S}_j(M) \) and \( \text{pd}_R \mathcal{S}_j(M) < \infty \) (see Theorem \ref{J1}).  

In Section \ref{section5}, we derive formulas for the Betti numbers of \( \mathcal{S}_j(M) \) in terms of those of \( M \) (see Proposition \ref{bettiprop} and Corollary \ref{cor 5.2}). As an application for linear-type modules, we give a formula for the Betti numbers of the \( j \)-th graded component of the Rees algebra \( \mathcal{R}_j(M) \). Furthermore, when \( M = I \), we provide an explicit formula for the Betti numbers of \( I^j \) in terms of those of \( I \) (see Corollary \ref{cor5.9}).  

Finally, in Section \ref{section6}, we establish upper and lower bounds for the Betti numbers of \( \mathcal{S}_j(M) \) (see Theorem \ref{Prop 6} and Corollary \ref{Bound1}). As elementary applications, we address special cases of the Buchsbaum-Eisenbud-Horrocks conjecture  (see Proposition \ref{prop 8}).  

The paper includes examples computed using Macaulay2 \cite{Macaulay}. 

\section{The Symmetric Algebra of a module finitely generated}\label{section2}

Let $R$ be a ring and $M$ an $R$-module, the \textit{symmetric algebra} of $M$, which we  denote by $\mathcal{S}_R(M)$, is the quotient of the tensor algebra $T_R(M)$ by the ideal
$<x\otimes y-y\otimes x;\, x, y\in M>$, that is,
$$\mathcal{S}_R(M):=\frac{T_R(M)}{<x\otimes y-y\otimes x; x,y\in M>}.$$
Since that $\mathcal{S}_R(M)$ is an graded $R$-algebra, i. e.,  $\mathcal{S}_R(M)=\underset{j\geq 0}{\oplus}\mathcal{S}_j(M)$, where $\mathcal{S}_0(M)=R$ and $\mathcal{S}_1(M)=M$. The $j$th component $\mathcal{S}_j(M)$ is called 
the \textit{$j$th symmetric power} of $M$.

The algebra $\mathcal{S}_R(M)$ can also be seen as an $R$-algebra together with a homomorphism of $R$-modules $\pi: M\longrightarrow \mathcal{S}_R( M)$ satisfying the following universal property: given an $R$-algebra $B$ and a homomorphism of $R$-modules $\psi: M\longrightarrow B$, there is a unique homomorphism of $R$-algebras $ \rho: \mathcal{S}_R(M)\longrightarrow B$ such that $\rho\circ\pi=\psi$, that is, such that the following diagram is commutative 
	$$
		\xymatrix{
			M\ar[r]^{\psi} \ar[d]_{\pi} & B \\
			\mathit{S}_R(M). \ar[ru]_{\rho}}
		$$

It is well known that, if $M$ is an $R$-module free of rank $n$,
$\mathcal{S}_R(M)$ is the ring of polynomials $R[T_1,\ldots, T_n].$

More generally, when $M$ is given by a presentation
$$R^m\overset{\phi}{\longrightarrow} R^l \overset{\pi}{\longrightarrow} M\longrightarrow 0, \\\ \phi=(a_{ij})_{l\times m},$$
its symmetric algebra is the quotient of the ring of polynomials $R[T_1,\ldots, T_l]$ by the ideal $Q$ generated by the $1$-forms $f_i= a_{1i}T_1+\cdots+a_{li}T_l, i = 1, 2,\ldots, m.$ Therefore, if $M$ is an $R$-module finitely generated over a Noetherian ring $R$ then $\mathcal{S}_R(M)$ is also finitely generated as $R$-module and the $j$th symmetric power $\mathcal{S}_j(M)$ is  also finitely generate as $R$-module. For more properties and results involving symmetric algebras (see \cite{Eisen}).

\section{Construction of complex symmetric power and results }\label{section3}

The purpose of this section is to construct a complex for each symmetric power of \(\mathcal{S}_R(M)\) associated with a free resolution of an \(R\)-module \(M\) (this complex was proposed by Tchernev in \cite{Tchernev}). Before proceeding with this construction, we make the following observations and considerations.

\begin{remark}\label{Char}
Let \(R\) be a local Noetherian ring and \(M\) a finitely generated \(R\)-module. Let \(F_\bullet\) be a finite free (minimal) resolution of \(M\) over \(R\). We denote the length of the complex \(F_\bullet\) by \(p\). Weyman, in \cite{Weyman}, initially studied the complexes of symmetric powers of an \(R\)-module \(M\) associated with the (minimal ) resolution \(F_\bullet\), denoted by \(\mathcal{S}_j(F_\bullet)\), and claimed that \(\mathcal{S}_j(F_\bullet)\) is a complex (or minimal free resolution). However, in \cite{15, 16}, it was later shown that the complexes \(\mathcal{S}_j(F_\bullet)\) defined in \cite{Weyman} are not, in fact, complexes. When verifying the condition \(d^2 = 0\), it is observed that there is a part of the differential that depends on two consecutive maps in \(F_\bullet\), and this part of the composition is not always zero.

Therefore, due to this problem, to ensure that \(\mathcal{S}_j(F_\bullet)\) is a complex, there are some ways to address the difficulty. One of them is to make assumptions about the characteristic of the ring \(R\), that is, to assume that the integers \(m\) satisfying \(2 \leq m \leq jp\) are invertible in \(R\). For example, \(R\) could be a \(K\)-algebra, where \(k\) is a field of sufficiently large characteristic. For more details on the properties of such complexes, see, for example, \cite{15, 16}.
\end{remark}

Due to Remark \ref{Char}, throughout this paper, our ring \( R \) will always satisfy assumptions about its characteristic. This is done to ensure that \( \mathcal{S}_j(F_\bullet) \) is a complex over \( R \).

Let $F$ be a finitely generated free $R$-module. The \textit{rank of $F$}, $\rank{F}$, is equal to the cardinality of the basis of $F$. If $F$ is free of rank $n$, then $F\cong R^n$.

Before introducing the complex, we need to define the $j$th \textit{divided power} of a free $R$-module $F$, denoted by $D_j(F)$ for $j\geq 0$ (for more details about see \cite{Buchs}).

\begin{definition}[Divided Power]
Let $F$ be a free $R$-module of finite rank, and $j\geq 0$ an integer nonnegative. The  \textit{$j$th divided power $D_j(F)$} is defined as the set of symmetric tensors in $T^j(F):=F^{\otimes j}$, that is,
$$D_j(F):=\{\omega\in T^j(F)|\hspace{0.1cm}\sigma(\omega)=\omega\hspace{0.1cm} \textrm{for all} \hspace{0.1cm}\sigma \in \mathfrak{S}_j\},$$
where $\mathfrak{S}_j$ is set the permutation of order $j.$ 
\end{definition}
Note that $D_0(F)=R$ and $D_1(F)=F$ and $D_j(F)$ is a free $R$-module. In fact, suppose $F$ be a finite free $R$-module generated by $f_1,\hspace{0.1cm}f_2, \ldots,\hspace{0.1cm}f_l$. To get a basis for $D_j(F)$ we first consider the orbits
$$\mathcal{O}_{a_1,a_2,\ldots,a_l}:=\mathfrak{S}_j\cdot f_1^{\otimes a_1}\otimes f_2^{\otimes a_2}\otimes\cdots\otimes f_l^{\otimes a_l}$$
and for $a_1+a_2+\cdots+a_l=j$, consider the divided power monomials
$$f_1^{(a_1)}\cdots f_l^{(a_l)}:=\underset{\omega \in \mathcal{O}_{a_1,\hspace{0.1cm}a_2,\ldots,\hspace{0.1cm}a_l}}{\sum} \omega$$
they form a basis for $D_j(F)$. i.e.,
\begin{equation}\label{D(F)}
     D_j(F)=\left\langle\{\underset{i}{\prod} f_i^{(a_i)}| \sum a_i=j\}\right\rangle
\end{equation}
is a free $R$-module.

\begin{lemma}\label{remarkD}
Let $F$ be a free $R$-module and $j$ be a nonnegative integer. If $\rank{F}= l$, then  $D_j(F)$ has rank
	$\displaystyle\binom{j+l-1}{l-1}$.
\end{lemma}
\begin{proof} By Equation \ref{D(F)}, the minimum number of generators of $D_j(F)$ can be seen to be the number of solutions, distinct, with nonnegative integers, of  the equation
$$a_1+a_2+\cdots+a_{l}=j$$ 
where $l=\rank{F}$. So, the number of solutions is $\displaystyle\test{l+j-1}{l-1}$. Thus, by \ref{D(F)}, we obtain the result.
\end{proof}

\subsection{Construction of $\mathcal{S}_j{\bf F}_{\bullet}$}\label{functor}

Let $R$ be a local Noetherian ring and $M$ be a finitely generated $R$-module. Assume that 
$${\bf F}_{\bullet}: 0\longrightarrow F_p\overset{\phi_p}{\longrightarrow} F_{p-1}\overset{\phi_{p-1}}{\longrightarrow}\cdots\overset{\phi_1}{\longrightarrow} F_0$$
is a finite free resolution of $M$ where $p:=\pd_RM$. Let $a_0,\, a_1,\ldots,\, a_p$ be a sequence of natural numbers. We define the functors
\begin{equation*}
 S(a_0,\ldots, a_p:\,{\bf F}_{\bullet}):=\left\{\begin{array}{rc} D_{a_0}F_0\otimes\overset{a_1}{\Lambda}F_1\otimes D_{a_2}F_2\otimes\cdots\otimes\overset{a_{p-1}}{\Lambda}F_{p-1}\otimes D_{a_p}F_p,&\mbox{for} \quad p\quad \mbox{even},
\\
D_{a_0}F_0\otimes\overset{a_1}{\Lambda}F_1\otimes D_{a_2}F_2\otimes\cdots\otimes D_{a_{p-1}}F_{p-1}\otimes\overset{a_p}{\Lambda}F_p, &\mbox{for}\quad p \quad \mbox{odd }.
\end{array}\right.
\end{equation*}
and the differential maps are as follows:
$$d^i:S(a_0,\ldots, a_p:{\bf F}_{\bullet})\rightarrow S(b_0,\ldots,\,b_p: {\bf F}_{\bullet})$$
which is defined to be zero when $(b_0,\ldots,b_p)\neq (a_0,\ldots,\hspace{0.1cm}a_i+1,\hspace{0.1cm}a_{i+1}-1,\ldots, a_p)$ for all $i$,  and in the case $(b_0,\dots, b_p)=(a_0,\dots, a_i+1,\hspace{0.1cm}a_{i+1}-1, \ldots,\hspace{0.1cm} a_p)$, 
\begin{equation}\label{differential1} 
\hspace{0.1cm} d^i=\left\{\begin{array}{rc}
 \pm 1\otimes\cdots \otimes 1\otimes A_{a_{i+1},\,a_i}\phi_{i+1}\otimes 1\cdots,&\mbox{for}\quad i\quad \mbox{odd} 
,\\
\pm 1\otimes\cdots \otimes 1\otimes B_{a_{i+1},\,a_i}\phi_{i+1}\otimes 1\cdots, &\mbox{for}\quad i \quad \mbox{even },
\end{array}\right. 
\end{equation}
\noindent where $\pm$ denotes $(-1)^{\sigma};$  $\sigma =a_0+2a_1+\cdots+(i+1)a_i$,\,\,  $A_{a_{i+1},a_i}\phi_{i+1}$ and $B_{a_{i+1},a_i}\phi_{i+1}$ are homomorphisms defined as follows: suppose that $f_1,\hspace{0.1cm}f_2,\ldots, f_r$ and $g_1,\hspace{0.1cm}g_2,\ldots, g_s$ form a basis for $F_{i+1} $ and $F_i$ respectively. Let
$$ A_{a_{i+1},a_i}\phi_{i+1}: D_{a_{i+1}}F_{i+1}\otimes \overset{a_i}{\Lambda} F_{i}\longrightarrow D_{a_{i+1}-1}F_{i+1}\otimes
\overset{a_i+1}{\Lambda} F_{i}$$
$$A_{a_{i+1},a_i}\phi_{i+1}(f_1^{(a_{i+1})_1}\cdots f_r^{(a_{i+1})_r}\otimes v)=\sum_{l=1}^r f_1^{(a_{i+1})_1}\cdots f_l^{(a_{i+1})_{l-1}}\cdots f_r^{(a_{i+1})_r}\otimes \phi_{i+1}(f_l)\wedge v$$
and
$$B_{a_{i+1},a_i}\phi_{i+1}:\overset{a_{i+1}}{\Lambda}F_{i+1}\otimes D_{a_i}F_i\longrightarrow \overset{a_{i+1}-1}{\Lambda}F_{i+1}\otimes D_{a_i+1}F_i$$
\small{$$B_{a_{i+1},i}\phi_{i+1}(f_{(a_{i+1})_1}\wedge \cdots \wedge f_{(a_{i+1})_s}\otimes w)=\displaystyle\sum_{l=1}^r (-1)^{l}f_{(a_{i+1})_1}\wedge \cdots \wedge \widehat{f_{(a_{i+1})_l}}\wedge \cdots \wedge f_{i_s}\otimes \phi_{i+1} (f_{(a_{i+1})_l})\cup w,$$}
\normalsize{where $g_i\cup g_1^{(i_1)}\cdots g_s^{(i_s)}= g_1^{(i_1)}\cdots g_i^{(i_i+1)}\cdots g_s^{(i_s)}$. Here 
$\widehat{f}_{i_l}$ means that $f_{i_l}$ is omitted.}

Thus, we define  $$(S_j{\bf F}_{\bullet})_t:=\underset{\sum ia_i=t}{\underset{\sum a_i=j}{\underset{(a_0, \dots, a_p)}{\oplus}}} S(a_0,\ldots, a_p:\,{\bf F}_{\bullet})\hspace{0.1cm}\,\,\,\mbox{for all }\,\,t\geq 0$$ 
and the differentials $d_{t}$ are given by
$$d_t:(S_j{\bf F}_{\bullet})_{t}\longrightarrow (S_j{\bf F}_{\bullet})_{t-1}\,\,\,\mbox{where }\,\,d_t:=(d_t^{j_1},\hspace{0.1cm}d_t^{j_2},\ldots,\hspace{0.1cm}d_t^{jr}),\hspace{0.1cm} \mbox{for all } \hspace{0.1cm} t\geq 1.$$
Thus, we get
\begin{equation}\label{equation2}
\mathcal{S}_j{\bf F}_{\bullet}:\cdots\longrightarrow (S_j{\bf F}_{\bullet})_{t+1}\overset{d_{t+1}}{\longrightarrow} (S_j{\bf F}_{\bullet})_t\overset{d_t}{\longrightarrow} \cdots\overset{d_2}{\longrightarrow}(\mathcal{S}_j{\bf F}_{\bullet})_1\overset{d_1}{\longrightarrow}(\mathcal{S}_j{\bf F}_{\bullet})_0.
\end{equation}
Note that, the notation $d_t^{j_r}$ is to indicate the differential $d^{j_r}$ on the $t$-th level of $\mathcal{S}_j{\bf F}_{\bullet}$ with $r$ the $r$-th solution of equation system
\begin{equation}\label{eqsoma}
    \sum_{i=0}^p ia_i=t;\,\,\,\,\,\,\sum_{i=0}^p a_i=j.
\end{equation}

\begin{remark}\label{construction} Summarizing the construction above, each component of $\mathcal{S}_j{\bf F}_{\bullet}$ is given by:
    $$(\mathcal{S}_j{\bf F}_{\bullet})_t =\begin{cases} \underset{\sum ia_i=t}{\underset{\sum a_i=j}{\underset{(a_0,\ldots, a_p)}{\oplus}}}D_{a_0}F_0\otimes\overset{a_1}{\Lambda}F_1\otimes D_{a_2}F_2\otimes\overset{a_3}{\Lambda}F_3\otimes\cdots\otimes D_{a_p}(F_p), &\mbox{if}\hspace{0.1cm}p\hspace{0.1cm} \mbox{is even},\\
\underset{\sum ia_i=t}{\underset{\sum a_i=j}{\underset{(a_0,\ldots, a_p)}{\oplus}}}D_{a_0}F_0\otimes\overset{a_1}{\Lambda}F_1\otimes D_{a_2}F_2\otimes\overset{a_3}{\Lambda}F_3\otimes\cdots\otimes\overset{a_p}{\Lambda}F_p, &\mbox{if}\hspace{0.1cm}p\hspace{0.1cm} \mbox{is odd}.
\end{cases}$$
\end{remark}
\begin{remark}(\cite[1.9]{Tchernev}).\label{remarklenght}
    Each $\mathcal{S}_j{\bf F}_{\bullet}$ in \ref{equation2} is a  bounded complex  as $R$-module and its length is given by
\begin{equation*}
\lambda(\mathcal{S}_j{\bf F}_{\bullet})=\begin{cases}
jp,&\mbox{for}\quad p \quad \mbox{even},\\

j(p-1)+\mbox{min}\{\rank{F_p},\hspace{0.1cm}j\}, &\mbox{for}\quad p \quad \mbox{odd}
\end{cases}
\end{equation*}
where   $\lambda(-)$ denotes $\lambda(\textbf{F}_{\bullet}):=\sup\{i|\, F_i\neq 0\}$ for some complex $\textbf{F}_{\bullet}$.
\end{remark}

\begin{remark}\label{pd1}
Suppose that $p=1$, then for each $j\geq 2$ the solutions of the the system of equations \ref{eqsoma} are given by $S=\{(a_0,\hspace{0.1cm}a_1)=(j-t,\, t)|\, t=0,\,1,\ldots,\,\lambda(\mathcal{S}_j\textbf{F}_{\bullet})\}$. Thus
$$(\mathcal{S}_j\textbf{F}_{\bullet})_t=\underset{a_1=t}{{\underset{a_0+a_1=j}{\underset{(a_0,\, a_1)}{\oplus}}}}D_{a_0}F_0\otimes\overset{a_1}{\Lambda}F_1=D_{j-t}F_0\otimes\overset{t}{\Lambda}F_1;$$
for all $t=0,\hspace{0.1cm}1,\ldots,\hspace{0.1cm} \lambda(\mathcal{S}_j\textbf{F}_{\bullet}).$ Hence $\mathcal{S}_j\textbf{F}_{\bullet}$ is given by
\small{$$\hspace{0.1cm}
\mathcal{S}_j{\bf F}_{\bullet}: 0\longrightarrow D_{j-l}(F_0)\otimes \overset{l}{\Lambda} F_1\overset{d_l}{\longrightarrow} D_{j-l+1}(F_0)\otimes \overset{l-1}{\Lambda} F_1\overset{d_{l-1}}{\longrightarrow} \cdots \overset{d_2}{\longrightarrow} D_{j-1}(F_0)\otimes \overset{1}{\Lambda} F_1\overset{d_1}{\longrightarrow }D_j(F_0)\otimes \overset{0}{\Lambda} F_1,$$}
where $l=\lambda(\mathcal{S}_j\textbf{F}_{\bullet})=\min\{\rank{F_1},\hspace{0.1cm} j\}.$
\end{remark}

Next, we define some notation to state the theorem that we use throughout the paper.
\begin{itemize}
\item[(a)]Let $A$ be an $m\times n$ matrix over $R$ where $m,n\geq 0$. For $t=1, 2,\dots, \mbox{min}\{m,n\}$ we denote by $I_t(A)$ the ideal generated by the $t$-minors of $A$ (the determinants of $t\times t$ submatrices). For systematic reasons one sets $I_t(A)=R$ for $t\leq 0$ and $I_t(A)= 0$ for $t>\mbox{min}\{m,n\}$. If $\phi: F\longrightarrow G$ is a finite free $R$-modules homomorphism, then $\phi$ is given by a matrix $A$ with respect to bases of $F$ and $G$. Therefore we may put $I_t(\phi)= I_t(A).$

 \item[(b)] Let $I$ be a proper ideal of $R$, then $\grade{(I)}$ denotes the length of maximal regular sequence contained in $I$.
\end{itemize}

\begin{theorem}{(\cite[Theorem 2.1]{Tchernev})}\label{th}
Let 
$${\bf F}_{\bullet}: 0\longrightarrow F_p\overset{\phi_p}{\longrightarrow} F_{p-1}\overset{\phi_{p-1}}{\longrightarrow}\cdots\overset{\phi_1}{\longrightarrow} F_0$$ be a finite free resolution
with $\coker(\phi_1)=M$ and $r_i=\displaystyle\sum_{n=i}^p(-1)^{n-i}\rank(F_n)$. Then, $\mathcal{S}_j{\bf F}_{\bullet}$ is exact if and only if
\begin{itemize}
    \item [(a)] $\grade(I_{r_i}(\phi_i))\geq ji,$ for all $i$ even, where $1\leq i\leq p$;
    \item [(b)] $\grade(I_{r_i-t}(\phi_i))\geq j(i-1)+1+t$ for $t=0,\, 1,\,\ldots, \,j-1$, for all
    $i$ odd, where $1\leq i\leq p$;
    \item [(c)] $j!$ is invertible in $R$.
\end{itemize}
If $\mathcal{S}_j{\bf F}_{\bullet}$ is exact for each $j$, it is a finite free resolution of the symmetric power $\mathcal{S}_j(M)$ for each $j$.
\end{theorem}

\begin{remark}

Note that when $M$ has projective dimension $1$, the assumption that $j!$ is invertible in $R$ is not necessary, for example, see \cite[Remark 2.6 and 2.7]{Tchernev}. Furthermore, for modules where $\pd_R M\geq 2$, it is possible to find examples where the construction $\mathcal{S}_j\textbf{F}$ does not even produce a complex when we fail to assume (c). That is, this construction does not, in general, produce a complex (see \cite[Example 6.2]{Tchernev}).
\end{remark}

\section{Minimal free Resolution for symmetric power}\label{section4}
In the previous section, Theorem \ref{th} established a criterion for the symmetric power $\mathcal{S}_j(M)$ to possess a free resolution. Knowing this, the focus of this section is to establish some criteria that will now ensure that $\mathcal{S}_j{\bf F}_{\bullet}$ be a minimal free resolution.

Hereafter, for simplicity, we use the definition below when we are talking about the exactness of $\mathcal{S}_j\textbf{F}_{\bullet}$.
\begin{definition}[\cite{Molica}]\label{SW_j}
Let $M$ be an $R$-module and let 
$${\bf F}_{\bullet}: 0\longrightarrow F_p\overset{\phi_p}{\longrightarrow} F_{p-1}\overset{\phi_{p-1}}{\longrightarrow}\cdots\overset{\phi_1}{\longrightarrow} F_0 $$
be a finite  minimal free resolution of $M$, i.e., $\textrm{Im} (\phi_i)\subset\mathfrak{m}F_{i}$ for all $i\in \mathbb{N}$. For each $j\geq 2$ we say that $M$ satisfies the \textit{$(SW_j)$ condition} if $\mathcal{S}_j{\bf F}_{\bullet}$ is a finite free resolution of $\mathcal{S}_j(M).$ 
\end{definition}

\noindent We remark that, by Theorem \ref{th}, we say that $M$ satisfies $(SW_j)$ condition if and only if 
\begin{itemize}
    \item[(a)] $\grade(I_{r_i}(\phi_i))\geq ji$, for all $i$ even, where $1\leq i\leq p$;
\item [(b)] $\grade(I_{r_i-t}(\phi_i))\geq j(i-1)+1+t$ for $t=0,\, 1,\,\ldots, \,j-1$, for all
    $i$ odd, where $1\leq i\leq p$;
    \item [(c)] $j!$ is invertible in $R$.
    
\end{itemize}

The example below illustrates an ideal $I$ that satisfies  $(SW_2)$ condition.
\begin{example}\label{Exem 1}
Let $R=k[[x,y,z]]$ be a ring of formal power series over a field $k$ and let  $I=(yz^2, x^2z, x^3y^2)$ be an ideal of $R$. By MACAULAY2, a minimal free resolution for $I$ is given by 
$\textbf{F}_{\bullet}:0\longrightarrow R^2\overset{\phi_1}{\longrightarrow }R^3\longrightarrow I\longrightarrow 0,$
\noindent where  $\phi_1$ is given by matrix $3\times 2$
\begin{center}
$[\phi_1]=\left(\begin{array}{cc}
            -yz & -xz^2  \\
            x^2&  0\\
            0 & z
        \end{array}\right).$\end{center}
We get that $I_1(\phi_1)=(-yz, x^2, -xz^2, z)$ and $I_2(\phi_1)=(x^3z^2, -yz^2,x^2z)$. Then, one obtain that $\grade{(I_1(\phi_1))}=2$ and $\grade{(I_2(\phi_1))}=1$. Hence, the ideal $I$ satisfies $(SW_2)$ condition.
\end{example}

By Theorem \ref{th}, we know that if the complex $\mathcal{S}_j{\bf F}_{\bullet}$ is exact, then $\mathcal{S}_j{\bf F}_{\bullet}$ is a free resolution of $R$-module $\mathcal{S}_j(M)$. Naturally, we can ask, when is $\mathcal{S}_j{\bf F}_{\bullet}$ a minimal free resolution.   In order to answer this question, we first show the following lemma.

\begin{lemma}\label{J2}
	Let $(R, \mathfrak{m}, k)$ be a Noetherian local ring and	$d^i$ be a map of free $ R$-modules defined in \ref{differential1}. Suppose that ${\textbf F}_{\bullet}$ is a minimal free resolution of $M$, then  $d^i \otimes 1_{\mathit{k}}=0$.
\end{lemma}
\begin{proof}
Let $1_k \colon k \to k$ be the identity map defined by $1_k(\overline{y}) = \overline{y}$. By definition of $d^i$, if $(b_0, b_1, \ldots, b_p) \neq (a_0, a_1, \ldots, a_i+1, a_{i+1}-1, \ldots, a_p)$ for all $i$, then $d^i = 0$ implies $d^i \otimes 1_k = 0$, and the result follows immediately in this case.

Now suppose that $(b_0, b_1, \ldots, b_p) = (a_0, a_1, \ldots, a_i+1, a_{i+1}-1, \ldots, a_p)$ for some $i$. We therefore need to consider two cases based on the parity of $i$: when $i$ is odd and when $i$ is even. We present the detailed argument for $i$ odd, as the case for $i$ even follows similarly.

\textbf{Case$(1)$. $i$ is odd}. Let
$d^i=\pm 1\otimes\cdots \otimes 1\otimes A_{a_{i+1},a_i}\phi_{i+1}\otimes 1\otimes\cdots$, where
 $$A_{a_{i+1},a_i}\phi_{i+1}: D_{a_{i+1}}F_{i+1}\otimes \overset{a_i}{\Lambda} F_{i}\longrightarrow D_{a_{i+1}-1}F_{i+1}\otimes
\overset{a_i+1}{\Lambda} F_{i}$$
is a map defined by
{\small{$$A_{a_{i+1},a_i}\phi_{i+1}(f_1^{(a_{i+1})_1}\cdots f_r^{(a_{i+1})_r}\otimes v)=\displaystyle\sum_{l=1}^r f_1^{(a_{i+1})_1}\cdots f_l^{(a_{i+1})_{l-1}}\cdots
f_r^{(a_{i+1})_r}\otimes \phi_{i+1}(f_l)\wedge v.$$}}

Let $\overline{y}\in k $. By hypothesis ${\bf F}_{\bullet}$ is a minimal free resolution of $M$  which implies that $\phi_{i+1}(f_l)=xf$ for some $x\in \mathfrak{m}$ and $f\in F_i$. So by the linearity of the tensor product, we get
$$
\hspace{-2.3cm}d^i\otimes 1_{k}(f_1^{(a_{{i+1}})_1}\cdots f_r^{(a_{i+1})_r}\otimes v\otimes\overline{y})=
d^i(f_1^{(a_{i+1})_1}\cdots f_r^{(a_{i+1})_r}\otimes v)\otimes 1_{k}(\overline{y})$$
$$\hspace{-3cm}=\pm 1\otimes\cdots \otimes 1\otimes A_{a_{i+1},a_i}\phi_{i+1}(f_1^{(a_{i+1})_1}\cdots f_r^{(a_{i+1})_r}\otimes v)\otimes 1\otimes\cdots\otimes 1_{k}(\overline{y})$$
$$\hspace{-0.7cm}=\pm 1\otimes\cdots\otimes 1\otimes \displaystyle\sum_{l=1}^r f_1^{(a_{i+1})_1}\cdots f_l^{(a_{i+1})_{l-1}}\cdots f_r^{(a_{i+1})_r}\otimes \phi_{i+1}(f_l)\wedge v\otimes 1\otimes\cdots\otimes 1_{k}(\overline{y})$$
$$\hspace{-1.7cm}=\pm 1\otimes\cdots\otimes 1\otimes \displaystyle\sum_{l=1}^r f_1^{(a_{i+1})_1}\cdots f_l^{(a_{i+1})_{l-1}}\cdots f_r^{(a_{i+1})_r}\otimes xf\wedge v\otimes 1\otimes\cdots\otimes 1_{k}(\overline{y})$$
$$\hspace{-0.8cm}=\pm 1\otimes\cdots\otimes 1\otimes \displaystyle\sum_{l=1}^r f_1^{(a_{i+1})_1}\cdots f_l^{(a_{i+1})_{l-1}}\cdots f_r^{(a_{i+1})_r}\otimes f\wedge v\otimes 1\otimes\cdots\otimes 1_{k}(\overline{xy})=0.$$
\end{proof}

\begin{remark}
Note that Tchernev and Weyman, in \cite{1}, provide a functorial construction of complexes \(LF_\bullet\) for any polynomial functor of degree \(j\), including symmetric power. This construction is functorial and characteristic independent but lacks minimality, which means that if the complex \(F_\bullet\) is minimal, the complex \(LF_\bullet\) is not. However, if \(F_\bullet\) has length \(p\), then \(LF_\bullet\) has length \(jp\). Below, we show that, under certain additional hypotheses and assumptions, $\pd_R \mathcal{S}_j(M)<\infty$ can also be minimal.
\end{remark}

\begin{theorem}\label{J1}
Let $M$ be a finitely generated $R$-module with  $\pd_RM<\infty$. If $M$ satisfies the $(SW_j)$ condition and ${\bf F}_{\bullet}$
is a minimal free resolution of $M$, then $\mathcal{S}_j{\bf F}_{\bullet}$ is a minimal free resolution of $\mathcal{S}_j(M)$ and  $\pd_R \mathcal{S}_j(M)<\infty$.
\end{theorem}
\begin{proof}
Let ${\bf F}_{\bullet}:0\longrightarrow F_p\longrightarrow F_{p-1}\longrightarrow\cdots\longrightarrow F_1\longrightarrow F_0$ be a minimal free resolution of $M$, where $p:=\pd_RM$. Since $M$ satisfies the $(SW_j)$ condition, by Theorem \ref{th}, the complex $\mathcal{S}_j \textbf{F}_{\bullet}$ is a free resolution for $\mathcal{S}_j(M)$. To show that $\mathcal{S}_j \textbf{F}_{\bullet}$ is a minimal free resolution it is enough to show that $d_t\otimes 1_k=0$ for all $t\geq 1$, where $d_t$ is a map defined by $$ d_t:(\mathcal{S}_j{\bf F}_{\bullet})_t\longrightarrow (\mathcal{S}_{j}{\bf F}_{\bullet})_{t-1}, \,\,\,\, d_t =(d_t^{j_1}, d_t^{j_2},\ldots,d_t^{jr}).$$

Now, let 
$$f_*=(f^{j_1}_*,\hspace{0.1cm}f^{j_2}_*,\ldots,\hspace{0.1cm}f^{j_r}_*)\in \underset{\sum ia_i=t}{\underset{\sum a_i=j}{\underset{(a_0,\ldots, \hspace{0.1cm}a_p)}{\oplus}}} S(a_0,\ldots, a_p;{\bf F}_{\bullet})\hspace{0.1cm},$$
where each $f^{j_r}_* \in S(a_0^r,\hspace{0.1cm}a_1^r,\ldots,\hspace{0.1cm} a_p^r)$ where $(a_0^r,\hspace{0.1cm}a_1^r,\ldots,\hspace{0.1cm}a_p^r)$ is the $r$th non-negative integer solution of system $\sum a_i=j,\sum ia_i=j$. Thus, as ${\bf F}_{\bullet}$ is a minimal free resolution, by Lemma \ref{J2} we obtain that $d_t^{j_r}\otimes 1_{k}=0$. Therefore,
$$
\begin{array}{lll}
 d_t\otimes 1_{k}(f^*,\overline{y}) &=& d_t(f^*)\otimes 1_{k}(\overline{y})\\
&=&(d_t^{j_1}(f_*^{j_1}),\hspace{0.1cm}d_t^{j_2}(f_*^{j_2}),\ldots,\hspace{0.1cm} d_t^{j^r}(f_*^{j_r}))\otimes 1_{k}(\overline{y})\\
&=&(d_t^{j_1}(f_*^{j_1})\otimes 1_{k}(\overline{y}),\hspace{0.1cm}d_t^{j_2}(f_*^{j_2})\otimes 1_{k}(u),\ldots,\hspace{0.1cm} d_t^{j^r}(f_*^{j_r})\otimes 1_{k}(\overline{y}))\\
&=&0
\end{array}
$$
\noindent for all $\overline{y}\in k$. 

\noindent Now $\pd_R \mathcal{S}_j(M)<\infty$ follows by construction of complex $\mathcal{S}_j{\bf F}_{\bullet}$ (see \ref{equation2}).
\end{proof}

\begin{remark}
Notice that in Theorem \ref{J1}, the conclusion that \(\pd_R M < \infty\) implies \(\pd_R \mathcal{S}_j(M) < \infty\) for all \(j \geq 1\). This implication holds regardless of any assumptions on the characteristic of the ring or the hypothesis \((SW_j)\). More specifically, if \(\pd_R M = p\), then it follows that \(\pd_R \mathcal{S}_j(M) \leq pj\) (see \cite[Theorems 12.1 and 13.1]{1}).
\end{remark}

Next, we give some examples and consequences. For this, we first give the definition of Betti numbers.

\begin{definition}
Let $(R,\mathfrak{m},k)$ be a Noetherian local ring and $M$ a finitely generated $R$-module with minimal free resolution given by
$$\textbf{F}_{\bullet}: \cdots\overset{\phi_i}{\longrightarrow}F_i\overset{\phi_{i-1}}{\longrightarrow}F_{i-1}\longrightarrow\cdots \overset{\phi_2}{\longrightarrow}F_1\overset{\phi_1}{\longrightarrow}F_0 $$
the invariants $\beta_i^R(M):=\rank{F_i}=\dim_k(\Tor_i^R(k, M))$ are called the\textit{ $i$-th Betti number of $M$}.
\end{definition}
\begin{remark}
Note that,  $\mathcal{S}_0{\bf F}_{\bullet}= R$ and  $\mathcal{S}_1{\bf F}_{\bullet}= {\bf F}_{\bullet}.$ In particular, if ${\bf F}_{\bullet}$ is a minimal free resolution of $M$, then $\beta_i^R(\mathcal{S}_1(M))=\beta_i^R(M),$ for all $i=0,\dots, p.
$ For this reason, we will always be considering the symmetric powers $\mathcal{S}_j(M)$ and $\mathcal{S}_j\textbf{F}_{\bullet}$ with $j\geq 2$.
\end{remark}

\begin{corollary}\label{cor} Let $M$ be an $R$-module finitely generated with $\pd_RM<\infty$. If $M$ satisfies the $(SW_j)$ condition , then 
\begin{equation*}
\pd_R \mathcal{S}_j(M)= \begin{cases}
j\pd_R M,&\mbox{for}\quad \pd_RM \quad \mbox{even},\\

j(\pd_RM-1)+\mbox{\rm min}\{\beta_{\pd_RM}^R(M),\hspace{0.1cm}j\}, &\mbox{for}\quad \pd_RM \quad \mbox{odd}.
\end{cases}
\end{equation*}
\end{corollary}
\begin{proof}
Since $\pd_RM<\infty$, we can consider  ${\bf F}_{\bullet}$ to be
a minimal free resolution of $M$. As $M$ satisfies the $(SW_j)$ condition  (Definition \ref{SW_j}),  we get $\mathcal{S}_j{\bf F}_{\bullet}$ is exact. Now, by Theorem \ref{J1},   $\mathcal{S}_j{\bf F}_{\bullet}$ is a minimal free resolution of $\mathcal{S}_j(M)$. Thus, $\pd_R \mathcal{S}_j(M)=\lambda(\mathcal{S}_j{\bf F}_{\bullet})$. Therefore, by Remark \ref{remarklenght}, we obtain the result.
\end{proof}

\begin{remark}\label{rem}
Let $M$ a finitely generated $R$-module with $0<\pd_RM<\infty$. If $M$ satisfies $(SW_j)$ condition, by Theorem \ref{J1} and the Auslander-Buchsbaum formula, we get $\label{eqsym}\pd_R(\mathcal{S}_j(M))\leq \dim R$. Thus, by Corollary \ref{cor}, we have the following cases:
\begin{enumerate}
    \item[(a)] Case $\pd_RM$ is even:  $\displaystyle j\leq \frac{\dim R}{\pd_R M}.$
    \item[(b)] Case $\pd_RM$ is odd: If $\mbox{\rm min}\{\beta_{\pd_RM}^R(M),\hspace{0.1cm}j\}=\beta_{\pd_RM}^R(M)$, then $ j\leq \frac{\dim R-\beta_{\pd_RM}^R(M)}{\pd_R M-1}$,\\
for $\pd_R M\neq 1$ and for $\pd_R M=1$, we get $\mbox{\rm min}\{\beta_1^R(M),\hspace{0.1cm}j\}\leq \dim R.$ Now, if $\pd_R M\neq 1$ and  $\mbox{\rm min}\{\beta_1^R(M),\hspace{0.1cm}j\}=j$, then $j\leq \frac{\dim R}{\pd_RM}.$
\end{enumerate}
\end{remark}
This show that the complexes $\mathcal{S}_j \textbf{F}_{\bullet}$, over a local ring of dimension $d$, do not always produce minimal free resolutions of $\mathcal{S}_j M$  for all $j$ (see Example \ref{ex4.7}).
\begin{example}\label{ex4.7}
Let $R=k[[x,y,z,w]]$ be a ring of formal power series over a field $k$ of characteristic zero and let $I=(xw, xz, yw, yz)$ be an ideal of $R.$
Since $j=3>\displaystyle\frac{\dim R}{\pd_R I}$, the complex $\mathcal{S}_3\textbf{F}_{\bullet}$ (in \ref{equation2}) does not produce a minimal free resolution for $\mathcal{S}_3(I)$. Firstly, computing in MACAULAY2, the minimal free resolutions of $I$ and $\mathcal{S}_3(I)$ are given respectively by
$$\textbf{F}_{\bullet}:0\longrightarrow R\overset{\phi_2}{\longrightarrow} R^4\overset{\phi_1}{\longrightarrow} R^4 \longrightarrow I\longrightarrow 0$$ and 
$$\mathcal{S}_3(I)_{\bullet}:0\longrightarrow R^4\longrightarrow R^{16}\longrightarrow R^{33}\longrightarrow R^{40}\longrightarrow R^{20}\longrightarrow \mathcal{S}_3(I)\longrightarrow 0,$$ 
where
$$[\phi_1]=\left(\begin{array}{ccccc}
            -y& 0 &-w & 0  \\
            x& 0& 0 & -w\\
            0& -y& z& 0 \\
            0& x& 0& z
        \end{array}\right), \,\,\,
        [\phi_2]=\left(\begin{array}{c}
           w\\
           -z\\
           -y\\
           x
        \end{array}\right).$$
      
On the other hand, now using the complex $\mathcal{S}_3\textbf{F}_{\bullet}$ (in \ref{equation2}), we get 
{\small{$$\mathcal{S}_3\textbf{F}_{\bullet}: 0\longrightarrow D_3(R)\longrightarrow  R^4\otimes D_2(R)\longrightarrow \overset{2}{\Lambda} R^4\otimes R\oplus R^4\otimes D_2(R)\longrightarrow R^4\otimes R\otimes R^4\oplus\overset{3}{\Lambda} R^4\longrightarrow $$}}
$$\longrightarrow R^4\otimes \overset{2}{\Lambda}R^4\oplus R^4\otimes R\longrightarrow D_2(R^4)\otimes R^4\longrightarrow D_3(R^4)\longrightarrow 0.$$
Since  $\grade{(I_1(\phi_2))}=4<6$, by Theorem \ref{th}, the ideal $I$ does not satisfy the $(SW_3)$ condition. Thus $\mathcal{S}_3\textbf{F}_{\bullet}$ is not an exact complex. Therefore, $\mathcal{S}_3\textbf{F}_{\bullet}$ is not a minimal free resolution of $\mathcal{S}_3(I)$.
\end{example}
According to Remark \ref{rem}  and Theorem \ref{J1}, a natural question arises. Are there intervals of $j$ where $\mathcal{S}_j\textbf{F}_{\bullet}$ is minimal free resolution for $\mathcal{S}_j(M)$? The next example illustrates that this can happen.
\begin{example}(\cite[Example 1.3]{C Huneke}) 
The ``generic" ideal of projective dimension one is given by the ideal defined by the exact sequence
$$\textbf{F}_{\bullet}:0\longrightarrow R^n\overset{\phi_1}{\longrightarrow }R^{n+1}\longrightarrow I\longrightarrow 0$$ 
where $[\phi_1]=(x_{rs})$ is a generic $n$ by $n+1$ matrix over a field $k$, and let $R=k[x_{rs}]$ be the polynomial ring over a field $k$. Let $I$ be an ideal of $R$  generated by the $n$ by $n$ minors of $[\phi_1]$. So, by \cite[Corollary 4]{H.E}, we have
\begin{equation}\label{eight1}
\grade{(I_{t}(\phi_1))}=(n-t+1)(n+2-t),\,\,\, t=1,\dots,\hspace{0.1cm}n.
\end{equation}
Checking the condition $(b)$ of the Theorem \ref{th}, with $n=r_1$, and using the equality \ref{eight1}, we get 
$$\grade{(I_{r_1-t}(\phi_1))}\geq 1+t \,\mbox{for } t=0,\,\ 1,\,\,\ldots,\,  j-1.$$
Thus, the ideal $I$ satisfies the condition $(b)$. Now, if $\textbf{F}_{\bullet}$ is a minimal free resolution, by Theorem \ref{J1} and Remark \ref{rem} item $(b)$, $\mathcal{S}_j\textbf{F}_{\bullet}$ is a minimal free resolution for $\mathcal{S}_j(I)$ for $j=1,\ldots,n.$
\end{example}

Motivated by this example, we have the following corollaries.
\begin{corollary}\label{coeJ1}
Let $M$ be a finitely generated $R$-module with $\pd_RM=1$. If $\grade{(I_j(\phi_1))}\geq \beta^R_1(M)-j+1$ for $j=1,\ldots,\hspace{0.1cm} \beta^R_1(M)$, then  $\mathcal{S}_j\textbf{F}_{\bullet}$ is a minimal free resolution of $\mathcal{S}_j(M)$, for $j=2,\ldots,\hspace{0.1cm} \beta^R_1(M)$. 
\end{corollary}
\begin{proof} Since $\pd_R M=1$,
by \cite[Proposition 3]{Avramov}, the $(SW_j)$ condition is equivalent to $\grade({I_j(\phi_1)})\geq \beta_1^R(M)-j+1$  for $j=1,\ldots, \beta^R_1(M)$. Now, from Theorem \ref{J1}, we have $\mathcal{S}_j\textbf{F}_{\bullet}$ is a minimal free resolution for $\mathcal{S}_j(M)$ for $j=2,\ldots,\hspace{0.1cm} \beta^R_1(M)$.

\end{proof}
Note that, in the Corollary \ref{coeJ1}, the number of grades to be tested, for $\mathcal{S}_j\textbf{F}_{\bullet}$ to be an exact sequence, is much smaller than proposed in the $(SW_j)$ condition. So, this leads us to formulate the following natural question. What is the minimum number of conditions for $\mathcal{S}_j\textbf{F}_{\bullet}$ to be exact?
\begin{corollary}\label{cor 4.9}
Let $R$ be a local Cohen-Macaulay domain and $M$ an $R$-module with a free minimal resolution,
$\textbf{F}_{\bullet}:0\longrightarrow R^m\overset{\phi_1}{\longrightarrow }R^{n}\longrightarrow M\longrightarrow 0$.
Let $I_t(\phi_1)$ denote the ideal in $R$ generated by the $t\times t$ minors of $[\phi_1]$ and $\mu(M_{\mathfrak{p}},R_{\mathfrak{p}})\leq n-m+\height(\mathfrak{p})-1$ for all non-zero primes $\mathfrak{p}$ in $R$. Then, the complex $\mathcal{S}_j\textbf{F}_{\bullet}$ is a minimal free resolution for $\mathcal{S}_j(M)$ for each $1\leq j\leq m.$
\end{corollary}
\begin{proof}
Suppose that $\mu(M_{\mathfrak{p}}, R_{\mathfrak{p}})\leq n-m+\height(\mathfrak{p})-1$ for all non-zero primes $\mathfrak{p}$ in $R$. By \cite[Theorem 1.1]{Huneke}, we get 
$$\height(I_t(\phi_1))\geq m+2-t \hspace{0.1cm}\mbox{for}\hspace{0.1cm} 1\leq t\leq m.$$
Now, since $R$ is Cohen-Macaulay (\cite[Corollary 2.1.4]{Bruns}),
\begin{equation}\label{eq4.10}
\grade{(I_t(\phi_1))}=\height(I_t(\phi_1))\geq m+2-t \hspace{0.1cm}\mbox{for}\hspace{0.1cm} 1\leq t\leq m.
\end{equation}
Checking the condition $(b)$ of the Theorem \ref{th}, with $m=r_1$ and using the inequality \ref{eq4.10}, we get the following inequalities
  $$\grade{(I_{r_1-t}(\phi_1))}\geq t+1\,\, \mbox{for}\,\,  1\leq t\leq m.$$
  Thus, $M$ satisfies the condition $(b)$. Now, by Theorem \ref{J1} and Remark \ref{rem} item $(b)$, $\mathcal{S}_j\textbf{F}_{\bullet}$ is a minimal free resolution for $\mathcal{S}_j(M)$ for $1\leq j\leq n.$
\end{proof}
\begin{corollary}\label{cor4HB}
Let $R$ be a local Cohen-Macaulay domain and $I$ an ideal of $R$ with a minimal free resolution 
$\textbf{F}_{\bullet}:0\longrightarrow R^n\overset{\phi_1}{\longrightarrow }R^{n+1}\longrightarrow I\longrightarrow 0$. 
If $\mu(I_{\mathfrak{p}}, R_{\mathfrak{p}})\leq \height{(\mathfrak{p})}$
for every non-zero prime $\mathfrak{p}$ in $R$, then the complex $\mathcal{S}_j\textbf{F}_{\bullet}$ is a minimal free resolution for $\mathcal{S}_j(I)$ for each $1\leq j\leq n.$
\end{corollary}
\begin{proof}
Follow immediately from Corollary \ref{cor 4.9}.
\end{proof}
\section{Betti numbers of symmetric power of modules}\label{section5}
The purpose of this section is to provide a partial affirmative answer to the question posed in the introduction. First, we present the following proposition.

\begin{proposition}\label{bettiprop}
Let $M$ be a finitely generated $R$-module with  $\pd_R M=:p<\infty$. Suppose that $\mathcal{S}_j{\bf F}_{\bullet}$ is a free resolution to $\mathcal{S}_j(M)$, then $\mathcal{S}_j{\bf F}_{\bullet}$ is a minimal free resolution if and only if 
\begin{itemize}
    \item [(a)] for $p$ even,
{\footnotesize{$$
    \beta_t^R(\mathcal{S}_j(M))=\hspace{-0.2cm}
\underset{\sum ia_i=t}{\underset{\sum a_i=j}{\underset{(a_0,\ldots,\hspace{0.1cm} a_p)}{\sum}}}\test{\beta_0^R(M)+a_0-1}{a_0}\test{\beta_1^R(M)}{a_1}\test{\beta_2^R(M)+a_2-1}{a_2}\test{\beta_3^R(M)}{a_3}\cdots\test{\beta_p^R(M)+a_p-1}{a_p},$$}}
    \item [(b)] for $p$ odd, 
    {\footnotesize{$$\beta_t^R(\mathcal{S}_j(M))=\underset{\sum ia_i=t}{\underset{\sum a_i=j}{\underset{(a_0,\ldots,\hspace{0.1cm} a_p)}{\displaystyle\sum}}}\test{\beta_0^R(M)+a_0-1}{a_0}\test{\beta_1^R(M)}{a_1}\test{\beta_2^R(M)+a_2-1}{a_2}\test{\beta^R_3(M)}{a_3} \cdots\binom{\beta_p^R(M)}{a_p},$$}}
\end{itemize}
\noindent for all $t=0,\hspace{0.1cm}1,\dots,\hspace{0.1cm} l:=\pd_R \mathcal{S}_j(M).$
\end{proposition}
\begin{proof}
Suppose that $p$ is even; in the case $p$ odd, the proof follows similarly. Let\\
${\bf F}_{\bullet}: 0\longrightarrow F_p\overset{\phi_p}{\longrightarrow} F_{p-1}\overset{\phi_{p-1}}{\longrightarrow}\cdots\overset{\phi_1}{\longrightarrow} F_0$  be a minimal free resolution of $M$. From (\ref{equation2}) we obtain the following complex
$$\mathcal{S}_j{\bf F}_{\bullet}:0\longrightarrow (S_j{\bf F}_{\bullet})_l\overset{d_l}{\longrightarrow} (S_j{\bf F}_{\bullet})_{l-1}\overset{d_{l-1}}{\longrightarrow} \cdots\overset{d_2}{\longrightarrow}(\mathcal{S}_j{\bf F}_{\bullet})_1\overset{d_1}{\longrightarrow}(\mathcal{S}_j{\bf F}_{\bullet})_0$$ for each integer $j\geq 2$. Now, by  Remarks  \ref{construction} and  Lemma \ref{remarkD},  we get that
\footnotesize{$$\rank (\mathcal{S}_j{\bf F}_{\bullet})_t=\hspace{-0.3cm}
\underset{\sum ia_i=t}{\underset{\sum a_i=j}{\underset{(a_0,\ldots, a_p)}{\sum}}}\test{\beta_0^R(M)+a_0-1}{a_0}\test{\beta_1^R(M)}{a_1}\test{\beta_2^R(M)+a_2-1}{a_2}\test{\beta_3^R(M)}{a_3}\cdots\test{\beta_p^R(M)+a_p-1}{a_p}$$} \noindent for all\normalsize{ $t=0,\, 1,\ldots,\,l$}.

\normalsize{Finally, the proof follows from the following fact: the free resolution $\mathcal{S}_j{\bf F}_{\bullet}$ is a  minimal free resolution for $\mathcal{S}_j(M)$ if and only if $\dim_k(\Tor^R_t(k,\mathcal{S}_j(M)))=\rank (\mathcal{S}_j{\bf F}_{\bullet})_t$  $t=0,\ldots,\hspace{0.1cm}l.$ }
\end{proof}
From Theorem \ref{J1} and Proposition \ref{bettiprop}, we get the next corollary.
\begin{corollary}\label{cor 5.2}
Let $M$ be a finitely generated $R$-module with  $\pd_R M=:p<\infty$. If $M$ satisfies $(SW_j)$ condition, then
\begin{itemize}
    \item [(a)] for $p$ even,
{\footnotesize{$$
    \beta_t^R(\mathcal{S}_j(M))=\hspace{-0.2cm}
\underset{\sum ia_i=t}{\underset{\sum a_i=j}{\underset{(a_0,\ldots, a_p)}{\sum}}}\test{\beta_0^R(M)+a_0-1}{a_0}\test{\beta_1^R(M)}{a_1}\test{\beta_2^R(M)+a_2-1}{a_2}\test{\beta_3^R(M)}{a_3}\cdots\test{\beta_p^R(M)+a_p-1}{a_p},$$}}

    \item [(b)] for $p$ odd, 
    {\footnotesize{$$\beta_t^R(\mathcal{S}_j(M))=\underset{\sum ia_i=t}{\underset{\sum a_i=j}{\underset{(a_0,\ldots,\hspace{0.1cm} a_p)}{\displaystyle\sum}}}\test{\beta_0^R(M)+a_0-1}{a_0}\test{\beta_1^R(M)}{a_1}\test{\beta_2^R(M)+a_2-1}{a_2}\test{\beta^R_3(M)}{a_3} \cdots\binom{\beta_p^R(M)}{a_p},$$}}
\end{itemize}
\noindent for all $t=0,1,\dots, \hspace{0.1cm}\pd_R \mathcal{S}_j(M).$
\end{corollary}

\begin{corollary}\label{corosix}
Let $M$ be a finitely generated $R$-module with projective dimension $1$ such that $\grade{(I_j(\phi_1))}\geq \beta^R_1(M)-j+1$  for all $j=1,\ldots,\hspace{0.1cm} \beta^R_1(M)$. Then
$$
\beta^R_t (\mathcal{S}_j(M))= \test{\beta^R_0(M)+j-t-1}{j-t}\test{\beta^R_1(M)}{t},
\hspace{0.1cm}\mbox{for all}\hspace{0.1cm} t=0,\hspace{0.1cm}1,\ldots,\hspace{0.1cm} \pd_R\mathcal{S}_j(M).$$
\end{corollary}

\begin{proof} Let ${\bf F}_{\bullet}$ be a minimal free resolution of $M$. By the assumption and by Corollary \ref{coeJ1}, the complex 
\small{$$\hspace{0.1cm}
\mathcal{S}_j{\bf F}_{\bullet}:0\longrightarrow D_{j-l}(F_0)\otimes \overset{l}{\Lambda} F_1\longrightarrow D_{j-l+1}(F_0)\otimes \overset{l-1}{\Lambda} F_1\longrightarrow \cdots \longrightarrow D_{j-1}(F_0)\otimes \overset{1}{\Lambda} F_1\longrightarrow D_j(F_0)\otimes \overset{0}{\Lambda} F_1 $$}
\normalsize{is a minimal free resolution for $\mathcal{S}_j(M),$ with $l=\pd_R\mathcal{S}_j(M)$ for all $j=1,\ldots, \hspace{0.1cm}\beta^R_1(M)$. Now, by Proposition \ref{bettiprop}, one  obtain 
$$\beta_t^R(\mathcal{S}_j(M))=\test{\beta^R_0(M)+j-t-1}{j-t}\test{\beta^R_1(M)}{t}, \hspace{0.1cm} \textrm{for all}\hspace{0.1cm} t=0,\hspace{0.1cm}1,\ldots,\hspace{0.1cm} l,$$ as we wanted to demonstrate.}
\end{proof}

\begin{corollary}\label{corosixt}
Let $M$ be a finitely generated $R$-module with $\pd_R M=2$ satisfying $(SW_j)$ condition. Then,
\begin{itemize}
    \item[(a)] for $j\geq t$,
$$\beta_t^R(\mathcal{S}_j(M))=
\overset{\lfloor\frac{t}{2}\rfloor}{\underset{r=0}{\displaystyle\sum}}\displaystyle \test{\beta^R_2(M)+r-1}{r}\test{\beta^R_1(M)}{t-2r}\test{\beta^R_0(M)+j-t+r-1}{j-t+r},$$
\item[(b)] for $j<t$,
$$\beta_t^R(\mathcal{S}_j(M))=\overset{\textrm{min}\{j,\lfloor\frac{t}{2}\rfloor\}}{\underset{r=t-j}{\displaystyle\sum}}\displaystyle \test{\beta^R_2(M)+r-1}{r}\binom{\beta^R_1(M)}{t-2r}\test{\beta^R_0(M)+j-t+r-1}{j-t+r},$$
for all $t=0,\hspace{0.1cm}1,\ldots,\hspace{0.1cm} \pd_R\mathcal{S}_j(M).$
\end{itemize}
\end{corollary}
\begin{proof}
Observe that, by Proposition \ref{bettiprop}, it is enough to calculate the non-negative integer solutions of the system
\begin{equation}\label{pd2}
\begin{cases} a_1+2a_2=t\\
 a_0+a_1+a_2=j
\end{cases}
\end{equation}
For this, we consider the following cases:  $j\geq t$ and $j<t$.
\begin{itemize}
    \item[(a)] Case $j\geq t$.
    Doing some calculations, it is possible to show that in this case the solutions of (\ref{pd2}) are given by triple $$(a_0,a_1,a_2)=(k-t+a_2,t-2a_2,a_2)\hspace{0.2cm} \textrm{where} \hspace{0.2cm} 0\leq a_2\leq \lfloor \frac{t}{2}\rfloor.$$
Since $M$ satisfies $(SW_j)$ condition, by Proposition \ref{bettiprop}, we obtain that  
$$\beta_t^R(\mathcal{S}_j(M))=
\overset{\lfloor\frac{t}{2}\rfloor}{\underset{r=0}{\displaystyle\sum}}\displaystyle \test{\beta^R_2(M)+r-1}{r}\test{\beta^R_1(M)}{t-2r}\test{\beta^R_0(M)+j-t+r-1}{j-t+r}.$$
\item[(b)] Case $j<t$. In this case, the solutions of (\ref{pd2}) are given by triple
$$
    (a_0,a_1,a_2)=(j-t+a_2,t-a_2,a_2)\hspace{0.2cm} \textrm{where}\hspace{0.2cm} t-j\leq a_2\leq \textrm{min} \{j,\lfloor \frac{t}{2}\rfloor\}.
$$
Thus, by Proposition \ref{bettiprop}, we get the result desired.
\end{itemize}
\end{proof}
The next example shows that the condition $(SW_j)$ of the Corollary \ref{cor 5.2} cannot be removed.
\begin{example}\label{exe 5.4}
Let $R=k[[x_1, x_2, x_3, y_1, y_2]]$ be a ring of formal power series over a field $k$ of characteristic zero and let
$I=(x_1x_2y_1, x_1x_3y_1, x_2x_3y_1, x_1x_2y_2, x_1x_3y_2, x_2x_3y_2)$ be an ideal of $R.$ By MACAULAY2, we obtain a minimal free resolution of $I$ given by 
$$\textbf{F}_{\bullet}:0\longrightarrow R^2\overset{\phi_2}{\longrightarrow} R^7\overset{\phi_1}{\longrightarrow} R^6\longrightarrow I\longrightarrow 0$$
where $\phi_1$ and $\phi_2$ have a matrix representation given by
$$[\phi_1]=\left(\begin{array}{ccccccc}
            -x_3& 0 &0 & 0& -y_2 & 0 & 0  \\
            x_2& -x_2& 0& 0& 0& -y_2& 0\\
            0& x_1& 0& 0& 0& 0&-y_2\\
            0& 0& -x_3& 0& y_1& 0& 0\\
            0 & 0& -x_2& x_2& 0 &y_1 &0\\
            0 &0& 0 &x_1 &0 &0 &y_1
        \end{array}\right) ;\,\,\,\,
        [\phi_2]=\left(\begin{array}{cc}
           y_2 & y_2\\
           0 & y_2\\
           -y_1 & -y_1\\
           0 & -y_1\\
           -x_3 & -x_3\\
           x_2 & 0\\
           0 & x_1
        \end{array}\right).$$
    Notice that
    $\grade{(I_4(\phi_1))}>2$, $\grade{(I_5(\phi_1))}>1$ and  $\grade {(I_2(\phi_2))}=4$. So, the ideal $I$ satisfies the $(SW_2)$ condition.
    Now, by Theorem \ref{J1}, we get that $\mathcal{S}_2\textbf{F}_{\bullet}$ is a minimal free resolution of $\mathcal{S}_2(I)$ given by
$$\mathcal{S}_2\textbf{F}_{\bullet}:  D_2(R^2)\longrightarrow R^7\otimes R^2\longrightarrow \overset{2}{\Lambda}R^7\oplus R^6\otimes R^2\longrightarrow R^6\otimes R^7\longrightarrow D_2(F_0).$$
On the other hand, by MACAULAY2, we get the minimal free resolution of $\mathcal{S}_2(I)$
$$\mathcal{S}_2(I)_{\bullet}: 0 \longrightarrow  R^3\longrightarrow R^{14}\longrightarrow R^{33}\longrightarrow R^{42}\longrightarrow R^{21} \longrightarrow \mathcal{S}_2(I)\longrightarrow 0.$$
We can see that the ranks of these two minimal free resolutions are the same, i.e., the Betti numbers are the same ( see  Table \ref{TB 1}).

Similarly, we compute the minimal free resolution of $\mathcal{S}_3(I).$ In fact, the complex $\mathcal{S}_3\textbf{F}_{\bullet}$ is given 
$$\mathcal{S}_3\textbf{F}_{\bullet}: 0\longrightarrow D_3(R^2)\longrightarrow  R^7\otimes D_2(R^2)\longrightarrow \overset{2}{\Lambda} R^7\otimes R\oplus R^6\otimes D_2(R^2)\longrightarrow $$
$$  R^7\otimes R^2\otimes R^6\oplus\overset{3}{\Lambda} R^7\longrightarrow D_2(R^6)\otimes R^2\oplus R^6\otimes \overset{2}{\Lambda}R^7\longrightarrow D_2(R^6)\otimes R^7\longrightarrow D_3(R^6)\longrightarrow 0.$$
Note that, by Theorem \ref{th}, the complex $\mathcal{S}_3\textbf{F}_\bullet$ is not a free resolution for $\mathcal{S}_3(I),$ because $\grade{(I_2(\phi_2))}<6$ (does not satisfy the $(SW_3)$ condition). On the other hand, by MACAULAY2, we get a minimal free resolution of $\mathcal{S}_3(I)$
$$\mathcal{S}_3(I)_{\bullet}:0\longrightarrow R^4\longrightarrow R^{32}\longrightarrow R^{97}\longrightarrow R^{160}\longrightarrow R^{146}\longrightarrow R^{56}\longrightarrow \mathcal{S}_3(I)\longrightarrow 0.$$
We see that the ranks coming from $\mathcal{S}_3\textbf{F}_{\bullet}$ do not coincide with the Betti numbers of a minimal free resolution of $\mathcal{S}_3(I)$( see Table \ref{TB 1} ).

    \begin{table}[h]
    \caption{}\label{TB 1}
\vspace{0.2cm}
\begin{tabular}{l|r|r}\hline
			
			$t$-th Betti Number & $\mathit{S}_2(I)$ & $\mathit{S}_3(I)$ \\ 
			\hline                               
			$0$ &    $21$ & $56$\\
			$1$ &    $42$ & $127$\\
			$2$ &     $33$  & $147$ \\
			$3$ &     $14$  &  $119$\\
			$4$ &     $3$ &  $60$\\
		\hline	

\end{tabular}
\end{table}
\end{example}

The following example shows that the invertibility of 2! (condition (c)) is essential in Tchernev's Theorem, as symmetric algebras change dramatically when we vary the characteristic.

\begin{example}
Let $R=k[[x_1,x_2,x_3,y_1,y_2]]$ be a formal power series ring over a field $k$, and let $I= (x_1x_2^2+x_2, x_2+y_1y_2, x_3+x_1^2, y_1^4+x_3, y_2^2)$ be an ideal of $R$. Furthermore, let $M=\bigwedge^2 I$ be a finitely generated $R$-module, where $\bigwedge^2 I$ denotes the second exterior power of $I$.

Using Macaulay2, we obtain the minimal free resolution of $\mathcal{S}_2(M)$ when $\text{char}(k)=0$:
$$0\longrightarrow R^{132}\longrightarrow R^{674}\longrightarrow R^{1329}\longrightarrow R^{1220}\longrightarrow R^{488}\longrightarrow R^{55}\longrightarrow M\longrightarrow 0. $$

However, when $\text{char}(k)=2$, the minimal free resolution becomes:
$$0\longrightarrow R^{133}\longrightarrow R^{678}\longrightarrow R^{1336}\longrightarrow R^{1226}\longrightarrow R^{490}\longrightarrow R^{55}\longrightarrow M\longrightarrow 0.$$

This shows that the Betti numbers differ significantly, as displayed in Table \ref{TB 1}.

\begin{table}
    \caption{Betti number comparison}
\vspace{0.2cm}
\begin{tabular}{l|r|r|r }\hline
			$i$-th Betti number & $\text{char}(k)=0$ & $\text{char}(k)=2$ \\
			\hline                               
			0 & 55 & 55 \\
			1 & 488 & 490 \\
			2 & 1220 & 1226 \\
			3 & 1329 & 1336 \\
			4 & 674 & 678 \\
			5 & 132 & 133 \\
		\hline	
\end{tabular}
\end{table}
\end{example}

\subsection{Applications: Betti Numbers of Linear Type Modules}

In this section, we present applications of our previous results. For linear type ideals $I$, we derive formulas expressing the Betti numbers of the power ideal $I^j$ in terms of the ideal $I$.

Let $R$ be a ring with total ring of fractions $Q$. The \textit{torsion submodule} of $M$ with respect to $R$ is defined as the kernel of the canonical map 
$$M \longrightarrow M \otimes_R Q,$$
denoted by $\mathcal{T}_R(M)$. When $\mathcal{T}_R(M) = 0$, we say $M$ is a \textit{torsion-free} $R$-module. If $\mathcal{T}_R(M) = M$, then $M$ is called a \textit{torsion} module. We say $M$ has \textit{rank} $r$ if $Q \otimes_R M$ is free of rank $r$ over $Q$ for some $r \in \mathbb{N}$.

\begin{definition}(\cite{Linear})
The \textit{Rees algebra} $\mathcal{R}(M)$ is defined as the image of the $R$-algebra homomorphism
$$\mathcal{S}_R(\sigma): \mathcal{S}_R(M) \longrightarrow \mathcal{S}_R(F).$$
Note that $\mathcal{S}_R(F)$ is a polynomial ring and $\mathcal{R}_R(M)$ forms its subalgebra. Since the $R$-torsion submodule $\mathcal{T}_R(\mathcal{S}_R(M))$ equals the kernel of $\mathcal{S}_R(\sigma)$ (\cite{S.U.V}), we have $\mathcal{R}(M) \cong \mathcal{S}_R(M)/\mathcal{T}_R(\mathcal{S}_R(M))$, showing that $\mathcal{R}(M)$ is independent of the choice of $\sigma$. We say $M$ is of \textit{linear type} if $\mathcal{T}_R(\mathcal{S}_R(M)) = 0$, i.e., when $\mathcal{S}_R(M) \cong \mathcal{R}_R(M)$.

For the following results, $\mathcal{R}_j(M)$ denotes the $j$-th graded component of $\mathcal{R}(M)$.

When $I$ is an ideal of $R$, we say $I$ is of \textit{linear type} if $\mathcal{S}_R(I) \cong \mathcal{R}_R(I)$, which is equivalent to $\mathcal{S}_j(I) \cong I^j$ for all $j$, or alternatively, to $\mathcal{T}_R(\mathcal{S}_R(I)) = 0$.
\end{definition}

\begin{proposition}\label{linearp}
Let $M$ be a linear type module of $R$ with $\pd_R M=p<\infty$. If $M$ satisfies $(SW_j)$ condition, then
\begin{itemize}
    \item [(a)] for $p$ even,
{\footnotesize{$$
    \beta_t^R(\mathcal{R}_j(M))=\hspace{-0.2cm}
\underset{\sum ia_i=t}{\underset{\sum a_i=j}{\underset{(a_0,\ldots,\hspace{0.1cm} a_p)}{\sum}}}\test{\beta_0^R(M)+a_0-1}{a_0}\test{\beta_1^R(M)}{a_1}\test{\beta_2^R(M)+a_2-1}{a_2}\test{\beta_3^R(M)}{a_3}\cdots\test{\beta_p^R(M)+a_p-1}{a_p},$$}}
    \item [(b)] for $p$ odd, 
    {\footnotesize{$$\beta_t^R(\mathcal{R}_j(M))=\underset{\sum ia_i=t}{\underset{\sum a_i=j}{\underset{(a_0,\ldots,\hspace{0.1cm} a_p)}{\displaystyle\sum}}}\test{\beta_0^R(M)+a_0-1}{a_0}\test{\beta_1^R(M)}{a_1}\test{\beta_2^R(M)+a_2-1}{a_2}\test{\beta^R_3(I)}{a_3} \cdots\binom{\beta_p^R(M)}{a_p},$$}}
\end{itemize}
\noindent for all $t=0,\hspace{0.1cm}1,\dots,\hspace{0.1cm} \pd_R\mathcal{R}_j(M).$

\end{proposition}

\begin{proof}
Since $M$ is a linear type module, we get the isomorphism $\mathcal{R}_R(M)\cong \mathcal{S}_R(M)$ of graded algebra, which implies that $\mathcal{S}_j(M)\cong \mathcal{R}_j(M)$ for all $j$. Thus, by Corollary \ref{cor 5.2}, we obtain the result.
\end{proof}

\begin{remark}
The canonical epimorphism from $\mathcal{S}_R(I)$ to $\mathcal{R}_R(I)$ is a homogeneous homomorphism of $\mathbb{N}$-graded algebras. An ideal $I$ is said to be {\it  2-syzygetic} if $\mathcal{S}_2(I)\cong I^2$. Furthermore $I$ is said to be {\it $m$-syzygetic} if $\mathcal{S}_j(I)\cong I^j$ for $1\leq j\leq  m$. 
\end{remark}

\begin{corollary}\label{cor5.9}
Let $I$ be a linear type ideal of $R$ (or $m$-syzygetic) with $\pd_R I=p<\infty$. If $I$ satisfies $(SW_j)$ condition, then
\begin{itemize}
    \item [(a)] for $p$ even,
{\footnotesize{$$
    \beta_t^R(I^j)=\hspace{-0.2cm}
\underset{\sum ia_i=t}{\underset{\sum a_i=j}{\underset{(a_0,\ldots,\hspace{0.1cm} a_p)}{\sum}}}\test{\beta_0^R(I)+a_0-1}{a_0}\test{\beta_1^R(I)}{a_1}\test{\beta_2^R(I)+a_2-1}{a_2}\test{\beta_3^R(I)}{a_3}\cdots\test{\beta_p^R(I)+a_p-1}{a_p},$$}}
    \item [(b)] for $p$ odd, 
    {\footnotesize{$$\beta_t^R(I^j)=\underset{\sum ia_i=t}{\underset{\sum a_i=j}{\underset{(a_0,\ldots, a_p)}{\displaystyle\sum}}}\test{\beta_0^R(I)+a_0-1}{a_0}\test{\beta_1^R(I)}{a_1}\test{\beta_2^R(I)+a_2-1}{a_2}\test{\beta^R_3(I)}{a_3} \cdots\binom{\beta_p^R(I)}{a_p},$$}}
\end{itemize}
\noindent for all $t=0,\hspace{0.1cm}1,\dots, \hspace{0.1cm}\pd_RI^j.$

\end{corollary}

\begin{proposition}
Let $M$ be an  $R$-module of linear type with projective dimension $1$. Then
$$
\beta^R_t (\mathcal{R}_j(M))= \test{\beta^R_0(M)+j-t-1}{j-t}\test{\beta^R_1(M)}{t},
$$
for all $t=0,\hspace{0.1cm}1,\ldots,\hspace{0.1cm} \pd_R\mathcal{R}_j(M)$.
\end{proposition}
\begin{proof}
Since $\pd_RM=1$, we get a minimal free resolution $$ \displaystyle \textbf{F}_{\bullet}:0\to R^{\beta_1^R(M)}\overset{\phi}{\to} R^{\beta_0^R(M)}\to M\to 0,$$ such that  $\rank{M}=\beta_0^R(M)-\beta_1^R(M).$ Furthermore, since $M$ is of linear type, $\mathcal{T}_R(\mathcal{S}_R(M))=0$. Thus, by \cite[Proposition 3]{Avramov}, we obtain that $\grade{(I_j(\phi))}\geq \beta_1^R(M)-j+2$ for all $1\leq j\leq \beta_1^R(M).$ Therefore, by Corollary \ref{cor 5.2} follow the result.
\end{proof}

\begin{corollary}
Let $M$ be an $R$-module of linear type with projective dimension $1$. Then
$$
\mu (\mathcal{R}_j(M))= \test{\mu(M)+j-1}{j},
$$
for all $t=0,\hspace{0.1cm}1,\ldots,\hspace{0.1cm} \pd_R\mathcal{R}_j(M)$.
\end{corollary}

\section{Upper and lower bounds}\label{section6}
In this section, we establish upper and lower bounds of Betti numbers of the $j$-th symmetric power of $M$.
\begin{theorem}\label{Prop 6}Let $M$ be a finitely generated $R$-module with $\pd_R M=:p>1$. If $M$ satisfies $(SW_j)$ condition, then
\begin{enumerate}
    \item[(a)] for $p$ even,
    $$\beta_t^R(\mathcal{S}_j(M))\leq \test{\displaystyle\sum_{i=0}^p\beta_i^R(M)+j(\frac{p+2}{2})}{j},$$
    \item[(b)] for $p$ odd,
     $$\beta_t^R(\mathcal{S}_j(M))\leq \test{\displaystyle\sum_{i=0}^p\beta_i^R(M)+j(\frac{p+1}{2})}{j},$$
 \noindent   for all  $t=0,\hspace{0.1cm}1,\ldots,\hspace{0.1cm} \pd_R\mathcal{S}_j(M).$  
\end{enumerate}

\end{theorem}
\begin{proof}
Suppose that $p$ is even; in the case $p$ odd, the proof follows similarly. By hypothesis $M$ satisfies $(SW_j)$ condition. So, by Corollary \ref{cor 5.2}, we obtain the following equality
\begin{equation}\label{6.2}
\begin{array}{lll}
    \beta_t^R(\mathcal{S}_j(M))&=&\hspace{-0.2cm}
\underset{\sum ia_i=t}{\underset{\sum a_i=j}{\underset{(a_0,\ldots, a_p)}{\sum}}}\test{\beta_0^R(M)+a_0-1}{a_0}\test{\beta_1^R(M)}{a_1}\test{\beta_2^R(M)+a_2-1}{a_2}\test{\beta_3^R(M)}{a_3}\cdots\test{\beta_p^R(M)+a_p-1}{a_p}\\\\
&\leq & \underset{\sum a_i=j}{\underset{(a_0,\ldots,\hspace{0.1cm} a_p)}{\sum}}\test{\beta_0^R(M)+a_0-1}{a_0}\test{\beta_1^R(M)}{a_1}\test{\beta_2^R(M)+
a_2-1}{a_2}\test{\beta_3^R(M)}{a_3}\cdots\test{\beta_p^R(M)+a_p-1}{a_p},
\end{array}
\end{equation}
where $t=0,\hspace{0.1cm}1,\dots,\hspace{0.1cm} \pd_R \mathcal{S}_j(M).$ Since $a_i\leq j$, for all  even integer positive $i$ between $0$ and $p$, we get that
\begin{equation}\label{ineq 7}
    \test{\beta_i^R(M)+a_i-1}{a_i}\leq\test{\beta_i^R(M)+j}{a_i}. 
\end{equation}
Now, from inequalities \ref{6.2} and \ref{ineq 7} we obtain that 
$$
\begin{array}{lll}
\beta_t^R(\mathcal{S}_j(M))&\leq & \underset{\sum a_i=j}{\underset{(a_0,\ldots, a_p)}{\sum}}\test{\beta_0^R(M)+j}{a_0}\test{\beta_1^R(M)}{a_1}\test{\beta_2^R(M)+j}{a_2}\test{\beta_3^R(M)}{a_3}\cdots\test{\beta_p^R(M)+j}{a_p}\\\\
&\leq &\test{\displaystyle\sum_{i=0}^p\beta_i^R(M)+j(\frac{p+2}{2})}{j}
\end{array}$$
\noindent 
where the last inequality follows by the Generalized Vandermonde\textquotesingle s identity.
\end{proof}
As an immediate consequence of Theorem \ref{Prop 6} above, we have the following corollary.
\begin{corollary}[Bound Upp-Low]\label{BUP}
Let $M$ be a finitely generated $R$-module with $\pd_R M=1$. If $M$ satisfies $(SW_j)$ condition, then 
 $$\test{\beta_1^R(M)}{t}\leq\beta_t^R(\mathcal{S}_j(M))\leq \test{\displaystyle\sum_{i=0}^p\beta_i^R(M)+j}{j},$$
 \noindent for all  $t=0,\hspace{0.1cm}1,\ldots,\hspace{0.1cm} \pd_R\mathcal{S}_j(M).$
\end{corollary}

As established in Corollary~\ref{corosix}, for modules of projective dimension $1$, the $(SW_j)$ condition admits a reformulation when the grade of the ideals $I_j(\phi_1)$ satisfies a particular inequality for all $j = 1, \ldots, \beta_1^R(M)$. Under these conditions, we obtain the following corollary.
\begin{corollary}\label{Bound1}
  Let  $R$ be a local ring and $M$ be a finitely generated $R$-module  with $\rm{pd}_R(M)=1$ such that 
  $\grade{(I_j(\phi_1))}\geq \beta^R_1(M)-j+1,$ for all $j=1,\ldots, \beta^R_1(M)$. Then,
  $$\test{\beta_1^R(M)}{t}\leq\beta_t^R(\mathcal{S}_j(M))\leq \test{\displaystyle\sum_{i=0}^p\beta_i^R(M)+j}{j},$$
\noindent   for all  $t=0,\hspace{0.1cm}1,\ldots,\hspace{0.1cm} \pd_R\mathcal{S}_j(M).$  
\end{corollary}

\begin{proposition}\label{Bound2}
Let $M$ be a finitely generated $R$-module with $\pd_R M=2$. If $M$ satisfies $(SW_j)$ condition, then
  \begin{enumerate}
      \item[(a)] If $j\geq t$, then $ 
\beta_t^R(\mathcal{S}_j(M))\geq \displaystyle\test{\beta_1^R(M)}{t}
,\hspace{0.1cm} \textrm{for all} \hspace{0.1cm} t=0,\hspace{0.1cm}1,\ldots,\hspace{0.1cm} j.
$
\item[(b)] If $j < t$, then $ 
\beta_t^R(\mathcal{S}_j(M))\geq \displaystyle\test{\beta_1^R(M)}{2j-t}
,\hspace{0.1cm} \textrm{for all} \hspace{0.1cm}
$ $t=j+1,\dots,\hspace{0.1cm} \pd_R \mathcal{S}_j(M).$
  \end{enumerate}

  \end{proposition} 
  \begin{proof}
Since $M$ satisfies $(SW_j)$ condition, by Theorem \ref{J1}, $\mathcal{S}_j \textbf{F}_{\bullet}$ is a minimal free resolution for $\mathcal{S}_j(M)$. So, by Corollary \ref{corosixt}, we need consider two cases $j\geq t$ and $j<t.$ 
\begin{itemize}
    \item[(a)] For $j\geq t$ and for all $t=0,\dots,j$, we obtain
    $$
  \begin{array}{lll}
      \beta_t^R(\mathcal{S}_j(M)) &=& \overset{\lfloor\frac{t}{2}\rfloor}{\underset{r=0}{\displaystyle\sum}}\displaystyle \test{\beta^R_2(M)+r-1}{r}\test{\beta^R_1(M)}{t-2r}\test{\beta^R_0(M)+j-t+r-1}{j-t+r}  \\
       &\geq & \overset{\lfloor\frac{t}{2}\rfloor}{\underset{r=0}{\displaystyle\sum}}\displaystyle\test{\beta_1^R(M)}{t-2r}\\
       &\geq & \displaystyle\test{\beta_1^R(M)}{t}.
  \end{array}  
    $$
\item[(b)] Similarly, for $j<t$ and  for all $t=j+1,\dots,\hspace{0.1cm} \pd_R \mathcal{S}_j(M)$, we get 
$$
\begin{array}{lll}
\beta_t^R(\mathcal{S}_j(M))&=&\overset{\textrm{min}\{j,\lfloor\frac{t}{2}\rfloor\}}{\underset{r=t-j}{\displaystyle\sum}}\displaystyle \test{\beta^R_2(M)+r-1}{r}\binom{\beta^R_1(M)}{t-2r}\test{\beta^R_0(M)+j-t+r-1}{j-t+r}\\\\
&\geq & \displaystyle\test{\beta_1^R(M)}{2j-t}.
\end{array}
$$
\end{itemize}
  \end{proof}

\subsection{Applications: Buchsbaum-Eisenbud-Horrocks (BEH) and Total Rank Conjecture (TR) Conjectures}

To present our applications, we first recall the celebrated Buchsbaum--Eisenbud--Horrocks conjecture (BEH) and Total Rank Conjecture (TR). 

\noindent {\bf Conjecture (BEH).} Let $(R, \mathfrak{m}, k)$ be a $d$-dimensional Noetherian local ring, and let $M$ be a nonzero finitely generated $R$-module. If $M$ has finite length and finite projective dimension, then for all $i\geq 0$, the Betti numbers of $M$ over $R$ satisfy
$$\beta_i^{R}(M)\geq \binom{d}{i}.$$

This conjecture is known to hold for local rings of dimension $\leq 4$ (see \cite{AB}), but remains open in higher dimensions. Partial results have been established in various cases, as shown in \cite{EG}, \cite{Cha}, \cite{San}, and \cite{Chang}.

\noindent {\bf Conjecture (TR).} Let $(R, \mathfrak{m}, k)$ be a $d$-dimensional Noetherian local ring, and let $M$ be a nonzero finitely generated $R$-module. If $M$ has finite length and finite projective dimension, then 
$$\sum_{i\geq 0}\beta_i^{R}(M)\geq 2^d.$$

The (TR) conjecture was proved by Avramov and Buchweitz \cite{AB} for local rings of dimension 5 containing their residue field. More recently, Walker \cite{walker} established the conjecture for complete intersection rings of odd characteristic, with subsequent work by VandeBogert and Walker (\cite{Walker1}) extending this to characteristic two.

We now present some direct results about (BEH) and (TR) conjectures when the module $M$ has projective dimension 1, without assuming $M$ has finite length. Specifically, the following results show that the Betti numbers of the $j$-th symmetric power satisfy the inequality in (BEH) conjectures.

\begin{proposition}\label{prop 8}
Let $R$ be a local ring of dimension $d$ and $M$ be a finitely generated $R$-module with $\pd_R M=1$. If $M$ satisfies the conditions $(SW_j)$ and $\beta_1^R(M)\geq d$, then
\[ 
\beta_t^R(\mathcal{S}_j(M))\geq \displaystyle\test{d}{t}
,\hspace{0.1cm} \textrm{for all}\hspace{0.1cm} t=0,\hspace{0.1cm}1,\ldots, \hspace{0.1cm} \pd_R\mathcal{S}_j(M).
\] 
\end{proposition}
\begin{proof}
Since $\pd_R M=1$ and the $(SW_j)$ condition holds, by Corollary \ref{cor 5.2}, we obtain
\begin{equation}\label{ine 7.6}
\beta_t^R(\mathcal{S}_j(M))\geq\test{\beta_1^R(M)}{t},\hspace{0.1cm}\mbox{for all}\hspace{0.1cm}t=0,1,\ldots, \pd_R\mathcal{S}_j(M).
\end{equation}
 Now, as $\beta_1^R(M)\geq d$, we get $\displaystyle\test{\beta_1^R(M)}{t}\geq\test{d}{t}$. Therefore, by inequality \ref{ine 7.6}, we obtain the result.
\end{proof}

The next corollary easily follows from Propositions \ref{BUP} and \ref{prop 8}.

\begin{corollary}
Let  $R$ be a local ring of dimension $d$ and $I$ be an ideal of linear type with projective dimension $1$ such that $\beta_1^R(I)\geq d$. Then
$$
\beta^R_t (I^j)\geq \test{d}{t},
$$
for all $t=0,1,\ldots, \pd_RI^j$.
\end{corollary}

\begin{remark}
    Proposition \ref{prop 8} requires the condition $\beta_1^R(M) \geq d$ for the desired inequality. This motivates the study of modules whose first Betti number is bounded below by the ring's dimension. Fiber product rings provide an example satisfying this condition. 

Let $(S,\mathfrak{s},k)$ and $(T,\mathfrak{t},k)$ be local rings with surjective homomorphisms $S\twoheadrightarrow k\twoheadleftarrow T$. The fiber product 
$$S\times_k T = \{(s,t) \mid \pi_S(s) = \pi_T(t)\}$$
is a Noetherian local ring having maximal ideal $\mathfrak{s}\oplus\mathfrak{t}$, residue field $k$, and which embeds as a subring of $S\times T$ \cite[Lemma 1.2]{AAM}. The fiber product is called non-trivial precisely when both $S \neq k$ and $T \neq k$. The fiber product rings considered here are non-trivial.

\end{remark}

\begin{proposition}\label{fiber}
Let $S\times_kT$  be a $d$-dimensional local ring.  Let $M$ be a finitely generated $S$-module with $\pd_{S\times_kT} M=1$ and satisfying $(SW_j)$ condition. Then
$$
\beta^{S\times_kT}_t (\mathcal{S}_j(M))\geq \test{\beta_1^T(k)}{t},
$$
for all $t=0,\hspace{0.1cm}1,\ldots,\hspace{0.1cm} \pd_{S\times_kT}\mathcal{S}_j(M)$.
\end{proposition}
\begin{proof} Since $\pd_{S\times_kT} M=1$ and satisfying $(SW_j)$, by Corollary \ref{cor 5.2}, we get
\begin{equation}\label{ine 7.7}
\beta^{S\times_kT}_t (\mathcal{S}_j(M))\geq \test{\beta^{S\times_kT}_1(M)}{t}, \,\,\,\mbox{for all} \hspace{0.1cm}t=0,\hspace{0.1cm}1,\ldots,\hspace{0.1cm} \pd_{S\times_kT}\mathcal{S}_j(M).
\end{equation}
Now, by \cite[Theorem 1.8]{moore}, one obtain that $\beta_1^{R\times_kS}(M)=\beta_0^S(M)\beta_1^T(k)+\beta_1^S(M).$ Thus, $\beta_1^{R\times_kS}(M)\geq \beta_1^T(k).$ Therefore, by inequality \ref{ine 7.7}, we obtain the new very special results for the (BEH) and (TR) Conjectures.
\end{proof}

\begin{corollary}\label{fiberBET}
Let $S\times_kT$  be a $d$-dimensional local ring with $d:=\dim (T)\geq \dim (S)$ and $T$ is a regular local ring.  Let $M$ be a finitely generated $S$-module with $\pd_{S\times_kT} M=1$ and satisfying $(SW_j)$ condition. Then
$$
\beta^{S\times_kT}_t (\mathcal{S}_j(M))\geq \test{d}{t},
$$
for all $t=0,\hspace{0.1cm}1,\ldots,\hspace{0.1cm} \pd_{S\times_kT}\mathcal{S}_j(M)$.
\end{corollary}

\begin{corollary}
Let $S\times_kT$  be a $d$-dimensional local ring with $d:=\dim (T)\geq \dim (S)$ and $T$ is a regular local ring.  Let $M$ be a finitely generated $S$-module with $\pd_{S\times_kT} M=1$ and satisfying $(SW_j)$ condition. Then
$$
\sum\limits_{t=0}^d\beta^{S\times_kT}_t (\mathcal{S}_j(M))\geq 2^d.
$$
\end{corollary}

\begin{remark}
By Corollary \ref{fiberBET}, suppose there is an ideal of linear type $I$ in $S\times_kT$,  with $\pd_RI=1$ and $d:=\dim (T)\geq \dim (S)$ where $T$ is a regular local ring. Then,
$$
\beta^{S\times_kT}_t (I^j)\geq \test{d}{t},
$$
for all $t=0,\hspace{0.1cm}1,\ldots,\hspace{0.1cm} \pd_{S\times_kT}I^j$.
\end{remark}
\noindent {\bf Acknowledgments:}
The authors would like to thank the referee for the careful reading and for his/her comments and suggestions, which have improved the quality of the paper.

\end{document}